# Building accurate initial models using gain functions for waveform inversion in the Laplace domain

Wansoo Ha[1] and Changsoo Shin[2]

[1] Department of Energy Resources Engineering, Pukyong National University

[2] Department of Energy Resources Engineering, Seoul National University

Corresponding author: Changsoo Shin (css@model.snu.ac.kr)

**Abstract**

We suggest an initial model building technique using time gain functions in the Laplace domain. Applying the gain expressed as a power of time is equivalent to taking the partial derivative of the Laplace-domain wavefield with respect to a damping constant. We construct an objective function, which minimizes the logarithmic differences between the gained field data and the partial derivative of the modeled data with respect to the damping constant. We calculate the modeled wavefield, the partial derivative wavefield, and the gradient direction in the Laplace domain using the analytic Green's function starting from a constant velocity model. This is an efficient method to generate an accurate initial model for a following Laplace-domain inversion. Numerical examples using two marine field datasets confirm that a starting model updated once from a scratch using the gradient direction calculated with the proposed method can be successfully used for a subsequent Laplace-domain inversion.

**Introduction**

Full waveform inversion is a promising method to recover subsurface information (Tarantola, 1984). It generally minimizes the differences between the recorded data and the modeled data by using a local-gradient based optimization method for efficiency (Virieux and Operto, 2009). A critical limitation of the local-gradient method is its dependency on the initial model. If the starting model is not close to the global minimum, an inversion can fall into a local minimum. A robust objective function or an accurate starting model is required for a successful inversion (Virieux and Operto, 2009).

Many researchers have tried to solve this problem. Traveltime tomography (Brenders and

Pratt, 2007; Operto et al., 2006; Zelt et al., 2005), migration velocity analysis (Al-Yahya, 1989; Symes, 2008), or streotomograph (Billette and Lambare, 1998) can be used to generate initial models for full waveform inversions. Acquiring wide azimuth data (Pratt et al., 1996; Ravaut et al., 2004) or inverting data sequentially starting from low frequency can mitigate the local minima problem (Bunks et al., 1995). Changing the objective function can make the inversion robust to a specific problem such as noise (Amundsen, 1991; Crase et al., 1990; Guitton and Symes, 2003).

As a kind of full waveform inversion method, Laplace-domain inversion usually minimizes the logarithmic differences between the observed and modeled data in the Laplace domain (Shin and Cha, 2008). It can recover macro-velocity models starting from homogeneous models as demonstrated by synthetic and field data examples (Koo et al., 2011; Park et al., 2013; Shin and Cha, 2008). However, we can expect a Laplace-domain inversion yield a better result provided with a better initial model. We are trying to develop an efficient initial model building method similar to the Laplace-domain inversion method in this research.

Laplace-domain full waveform inversion naturally puts large weight on the early-arrival signal due to the damping in the Laplace transform (Shin and Cha, 2008). It is robust to the initial guess and results in large-scale subsurface background velocity models by sacrificing the late-arrival signal (Ha and Shin, 2013). However, the late-arrival signal also contains valuable information about the subsurface. Kwak et al. (2013) showed that the late-arrival can be used to enhance the results of a synthetic Laplace-Fourier domain full waveform inversion by applying time windows to the wavefield. Temporal gain was applied to the wavefield to generate time windows in the Laplace-Fourier domain (Kwak et al., 2013). Applying a gain function expressed as $t^n$ to the time domain wavefield is equivalent to calculating the partial derivative of the wavefield with respect to the damping constant in the Laplace domain (Kreyszig, 2011). However, numerical calculation of the derivative in the full waveform inversion involves recursive propagation of partial derivative wavefields and it makes the inversion computationally burdensome (Kwak et al., 2013).

In this study, we apply gain functions to the time domain wavefield and calculate the partial derivative wavefield in the Laplace domain using the analytic Green's function in the Laplace domain (Ha et al., 2011). This process is very cheap compared with the numerical calculation of the partial derivative wavefield (Kwak et al., 2013). However, we cannot iterate the process since we use the analytic solution starting from a homogeneous starting model. Instead, the model updated once by this method can be used as a better initial guess of the subsequent Laplace-domain full waveform inversion than a scratch model. We demonstrate the method using two marine field data examples.

**Effect of gain functions on the Laplace-domain wavefields**

We review the effect of the damping function in the Laplace transform and the gain function on the Laplace-domain wavefields. Laplace-domain full waveform inversion inherently utilizes a damped wavefield (Shin and Cha, 2008). The Laplace transform is defined as

$$\tilde{d}(s) = \int_0^\infty d(t) e^{-st} dt, \qquad (1)$$

where $d(t)$ is the time-domain wavefield, $\tilde{d}(s)$ is the Laplace-domain wavefield, and $s$ is a positive damping constant. The equation indicates that the Laplace-domain wavefield is equivalent to the integral of the damped wavefield damped by an exponential function. Figure 1a shows a shot gather from a marine survey at the Gulf of Mexico, and figure 1b shows a damped shot gather damped by $e^{-2t}$. Note that the late-arrival signal is disappearing even for a small damping constant of 2 $s^{-1}$. In this way, a Laplace-domain inversion uses the early-arrival signal to extract background subsurface velocity information. Therefore, the penetration depth of the Laplace domain inversion is shallow unless the maximum offset is large (Ha et al., 2012a).

However, the late-arrival signal also contains information about deeper structures. A gain function can be used to mitigate the elimination of the late-arrival signal in the damped wavefield. There are many kind of gain functions; however, we apply a power of time, $t^n$, as the gain function to exploit the function in a Laplace domain inversion. Figure 1c shows the shot gather with the gain function of $t^4$. By applying the damping function to the gained seismogram, we can obtain the shot gather in figure 1d. The gained and damped shot gather (Figure 1d) contains more late-arrival signal when compared with the original damped shot gather (Figure 1b).

Figure 2 shows the Laplace-domain wavefields with three different gain functions for the same damping constant. The amplitude of the large-offset signal increases as gain increases, because large-offset signal mainly contains the late-arrival signal. Therefore, we are putting more weight on the late-arrival signal as we increase the gain. We will use the gained seismogram to construct an initial model for a Laplace-domain full waveform inversion, which is better than homogeneous initial models used in many Laplace-domain researches (Koo et al., 2011; Park et al., 2013; Shin and Cha, 2008).

**Theory**

The Laplace domain wavefields can be obtained by Laplace transforming the time domain wavefields as

$$\tilde{u}(s) = \int_0^\infty u(t)e^{-st}\,dt,$$

$$\tilde{d}(s) = \int_0^\infty d(t)e^{-st}\,dt, \qquad (2)$$

where $u$ is the modeled wavefield, $d$ is the observed wavefield, and $s$ is a positive damping constant. Applying a gain function $t^n$ to the observed wavefield is equivalent to taking the partial derivative to the Laplace domain wavefield with respect to the damping constant (Kreyszig, 2011) as

$$\int_0^\infty d(t)t^n e^{-st}\,dt = (-1)^n \frac{\partial^n \tilde{d}(s)}{\partial s^n}. \qquad (3)$$

The logarithmic objective function of a Laplace-domain full waveform inversion (Shin and Cha, 2008) minimizes the logarithmic differences between the observed and modeld wavefield as

$$E(s) = \frac{1}{2}\sum_i^{N_s}\sum_j^{N_r}\left[\ln\frac{\tilde{u}_{ij}}{\tilde{d}_{ij}}\right]^2. \qquad (4)$$

We can build an objective function using the partial derivative wavefield in the Laplace domain. This objective function minimizes the logarithmic differences between the partial derivative of the observed and modeled wavefield as

$$E(s,n) = \frac{1}{2}\sum_i^{N_s}\sum_j^{N_r}\left[\ln\left(\frac{\partial^n \tilde{u}_{ij}}{\partial s^n}\bigg/\frac{\partial^n \tilde{d}_{ij}}{\partial s^n}\right)\right]^2. \qquad (5)$$

Theoretically, the objective function above also minimizes the differences between the gained wavefields of the observed and modeled data in the time domain. The gradient direction of the objective function using the partial derivative wavefield can be expressed as

$$\frac{\partial E(s,n)}{\partial m_k} = \sum_i^{N_s} \sum_j^{N_r} \frac{1}{\partial^n \tilde{u}_{ij} / \partial s^n} \left( \frac{\partial^{n+1} \tilde{u}_{ij}}{\partial m_k \partial s^n} \right) \ln \left( \frac{\partial^n \tilde{u}_{ij}}{\partial s^n} \Big/ \frac{\partial^n \tilde{d}_{ij}}{\partial s^n} \right), \qquad (6)$$

where, $m_k$ is the $k$ th model parameter, $N_s$ is the number of shots, and $N_r$ is the number of receivers. We use the analytic Green's function to calculate the gradient direction and update the subsurface model only once from a homogeneous model to obtain an initial model for a following Laplace-domain full waveform inversion. The Green's function of the acoustic wave equation in the Laplace domain can be expressed as

$$G = \frac{1}{4\pi} \left( \frac{e^{-st_1}}{r_1} - \frac{e^{-st_2}}{r_2} \right), \qquad (7)$$

where, $t_1 = r_1/c$, $t_2 = r_2/c$, and $r_1$ is the distance between a source and a receiver, $r_2$ is the distance between an imaginary source and the receiver, and $c$ is the velocity of the medium. We added the Lloyd mirror effect to consider the free surface boundary (Officer, 1958). Therefore, the observed wavefield in the Laplace domain can be expressed as

$$\tilde{u}_{ij} = \frac{f_m}{4\pi} \left( \frac{e^{-st_1}}{r_1} - \frac{e^{-st_2}}{r_2} \right), \qquad (8)$$

where, $f_m$ is the Laplace-domain source used to generate the modeled data. The source wavelet can be estimated using the observed and modeled data without applying the gain function (Shin and Cha, 2008). The partial derivative of the modeled wavefield with respect to the damping constant can be calculated as

$$\frac{\partial^n \tilde{u}_{ij}}{\partial s^n} = \frac{f_m}{4\pi} \frac{\partial^n}{\partial s^n} \left( \frac{e^{-st_1}}{r_1} - \frac{e^{-st_2}}{r_2} \right) = \frac{f_m}{4\pi} \left[ (-st_1)^n \frac{e^{-st_1}}{r_1} - (-st_2)^n \frac{e^{-st_2}}{r_2} \right]. \quad (9)$$

Calculation of the gradient direction (equation 6) also requires the partial derivative of the wavefield with respect to the model parameter. This can be calculated by taking the partial derivative of the wave equation in the Laplace domain. The acoustic Laplace domain wave equation (Shin and Cha, 2008) can be expressed as

$$\mathbf{S}\tilde{\mathbf{u}}_i = \left( \frac{s^2}{c^2} \mathbf{M} + \mathbf{K} \right) \tilde{\mathbf{u}}_i = \tilde{\mathbf{f}}_i, \quad (10)$$

where, $\mathbf{S}$ is the impedance matrix, $\mathbf{M}$ is the mass matrix, $\mathbf{K}$ is the stiffness matrix, and $\tilde{\mathbf{f}}_i$ is the Laplace-domain source vector. By applying the partial derivative to the wave equation with respect to the $k$ th model parameter, we can obtain

$$\frac{\partial \tilde{\mathbf{u}}_i}{\partial m_k} = \mathbf{S}^{-1} \left( -\frac{\partial \mathbf{S}}{\partial m_k} \tilde{\mathbf{u}}_i \right). \quad (11)$$

The partial derivative wavefield can be calculated using the analytic solution as

$$\frac{\partial \tilde{u}_{ij}}{\partial m_k} = \mathbf{S}^{-1} \left( \frac{2s^2}{c^3} \tilde{u}_{ik} \right) = \frac{f_m}{(4\pi)^2} \left( \frac{2s^2}{c^3} \right) \left( \frac{e^{-st_3}}{r_3} - \frac{e^{-st_4}}{r_4} \right) \left( \frac{e^{-st_1}}{r_1} - \frac{e^{-st_2}}{r_2} \right) = \frac{f_m s^2}{8\pi^2 c^3} \sum_{l=1}^{4} \frac{e^{-s\tau_l}}{\rho_l}, \quad (12)$$

where,

$$\tau_1 = t_1 + t_3, \ \tau_2 = t_2 + t_3, \ \tau_3 = t_1 + t_4, \ \tau_4 = t_2 + t_4,$$

$$\rho_1 = r_1 r_3, \ \rho_2 = -r_2 r_3, \ \rho_3 = -r_1 r_4, \ \rho_4 = r_2 r_4. \quad (13)$$

Therefore, the partial derivative wavefield with respect to both the subsurface parameter and the damping constant can be obtained by applying the partial derivative to the equation above with respect to the damping constant as

$$\frac{\partial^{n+1} \tilde{u}_{ij}}{\partial m_k \partial s^n} = \frac{f_m}{8\pi^2 c^3} \sum_{l=1}^{4} \left[ n(n-1)(-\tau_l)^{n-2} + 2ns(-\tau_l)^{n-1} + s^2 (-\tau_l)^n \right] \frac{e^{-s\tau_l}}{\rho_l}. \quad (14)$$

The gradient direction (eqaution 6) of the objective function can be obtained by multiplying the inverse of equation 9, equation 12 and the logarithmic differences between the partial derivative wavefields of the modeled and observed data. The final velocity update direction can be obtained by regularizing the gradient direction by the Hessian (Ha et al., 2012b; Pratt et al., 1998).

**Numerical examples**

We applied the gain method to a marine field dataset acquired at the Gulf of Mexico (Figure 1a). The dataset contains 399 shots each with 408 receivers. The maximum offset is 10,321 m and the minimum offset is 137 m. The recording time is 12 s and the sampling rate is 4 ms. The shot interval is 50 m and the receiver interval is 25 m.

Figure 3 show the gradient direction obtained using different power values from zero to four in the gain function. We calculated the gradient direction using a homogeneous starting model with the velocity of 3.5 km/s. The grid size used for the gradient calculation is 200 m, and we used every 4th shot gather for efficiency. Note that figure 3a is equivalent to the first gradient of the original Laplace-domain inversion since we applied no gain to the data. We can see shape of salt top and sedimentary layers below the water bottom. Since the gain function weights late-arrival signal, which contains information of the deeper structures, we can see the penetration depth of the gradient deepens as the power increases (Figure 1a to e). Figure 4a shows a weighted sum of the gradient directions shown in figure 1. We controled the weights to make each gradient contributes to the summed gradient equally. Figure 1b shows the homogeneous starting model with the water layer and figure 1c shows the updated model using the summed gradient from the homogeneous model. When we update the velocity model, we used the parabolic fitting to find the optimal step length (Press et al., 1992). Note that the calculation is cheap when compared with the numerical inversion methods because we used the analytic solution with the large grids.

We used the velocity models shown in figures 1b and c as the initial models for subsequent Laplace-domain full waveform inversions to examine the usefulness of the initial model generated by the gain method (Shin and Cha, 2008). The two inversions use exactly same setting except for the initial model. The grid size used in the inversion is 25 m. We interpolated the initial model by the gain method linearly to fit the model to the grids of the inversion. We inverted 11 damping constants simultaneously from 2 $s^{-1}$ to 12 $s^{-1}$. Figures 5 shows the two inversion results after 80 iterations. We can see that the artifacts below the water layer above the salt top is diminished in the velocity model obtained using the updated initial model. The error histories show that the initial model using the gain method significantly reduces the initial error and makes the inversion converge faster (Figure 6).

We applied the proposed method to a second field dataset. The dataset contains 1,156 shots with the interval of 37.5 m. Each shot gather contains 804 receivers with the interval of 12.5 m. The minimum offset is 165 m and the maximum offset is 10,202.5 m. The sampling rate is 4 ms and the recording time is 15 s. Figure 7 shows the gradient directions calculated using the analytic Green's function for the power of zero to four. We used every 4th shot gather with the grid size of 200 m to calculate the gradient directions. We can see the penetration depth of the gradient deepens as the power increases as the first example. It is hard to see salt structures at depth from the gradients;

however, the gradient direction recovers the sedimentary layer below the water layer. The weighted sum of each gradient (Figure 8a) is used to obtain the updated initial model from the homogeneous velocity model (Figure 8b) with the velocity of 3.5 km/s. Figure 8c shows the updated model using the parabolic fitting (Press et al., 1992). We used the homogeneous model and the updated model as the starting models for Laplace-domain full waveform inversions. We used 25 m grids and 11 damping constants ranging from 2 $s^{-1}$ to 12 $s^{-1}$. In this example, the inversion results using the two initial models are similar to each other (Figure 9). However, the error histories show that the initial model by the gain method makes the inversion start from a point closer to a minimum of the objective function than the homogeneous model (Figure 10). It reduces the amount of error and accelerates the convergence.

**Discussions**

The partial derivative wavefield with respect to the damping constant can be calculated numerically; however, the partial derivative wavefield of order *n* requires *n* more modeling for each shot, and it makes the algorithm computationally intensive even for 2D inversions of field data (Kwak et al., 2013).

We used the analytic Green's function to calculate the partial derivative wavefield used in the gradient calculation (equation 6). Since we used the analytic solution, we can obtain the first gradient only using a constant velocity model. The computational burden of this method is ignorable when compared with the numerical approach. We used the method to generate an accurate starting model for following Laplace-domain full waveform inversions. Field data examples showed that the initial model obtained using the gain method can be used as an accurate starting model for Laplace-domain inversion successfully.

We used the powers of time as the gain function in the time domain, which is equivalent to the partial derivative of the wavefield with respect to the damping constant in the Laplace domain. A common alternative gain function is the exponential function. An exponential function of time, $e^{rt}$, can be used as the gain function, with a positive $r$; however, it changes the damping constant due to the exponential damping in the Laplace transform as

$$\int_0^\infty d(t)e^{rt}e^{-st}\,dt = \int_0^\infty d(t)e^{-(s-r)t}\,dt = \tilde{d}(s-r). \tag{15}$$

Therefore, this is equivalent to decreasing the damping constants in a Laplace-domain inversion. A discussion about the range of the damping constant is given in (Ha et al., 2012a).

The power value we used in the gain function of numerical examples can be varied. However, large power value can cause instability in the Laplace transform. Figure 11a shows a trace from the first field data with the offset of 1,021 m. For a stable Laplace transform, the amplitude of the damped trace at the maximum recording time need to be small enough. We applied the exponential damping with the damping constant of 2 $s^{-1}$ to the original trace and two gained traces with the power of 4 and 8 (Figure 11b) in the gain function (equation 3). When the power is large, the amplitude of the last sample is not ignorable even with single precision calculation. Therefore, we need to limit the maximum power of the gain function for a stable Laplace transform of the recorded data. The maximum value depends on several parameters including the overall amplitude level of the data, the maximum recording time, a desired precision for the Laplace transform, and the maximum depth we want to recover (Ha et al., 2012a).

**Conclusions**

We proposed an efficient method to build a starting model for Laplace-domain full waveform inversions. The method used the partial derivative of the wavefield with respect to the damping constant in the Laplace domain, which is expressed as a gain in the time domain. The proposed objective function minimizes the differences between the partial derivative wavefields of the observed and modeled data. The partial derivative wavefield of the observed data was generated by Laplace-transforming the gained data with a power of time. The partial derivative wavefield of the modeled data was generated from the analytic Green's function. This is a straightforward process since we cannot iterate using the analytic solution. The resultant model can be used for following processes such as a Laplace-domain inversion. Laplace-domain inversions generally yield good results even with scratch initial models. However, the initial model calculated efficiently using the proposed method can enhance the results and convergence characteristics of following Laplace-domain inversions as shown by field data examples.

**References**

Al-Yahya, K., 1989, Velocity analysis by iterative profile migration: Geophysics, **54**, 718-729.

Amundsen, L., 1991, Comparison of the least-squares criterion and the Cauchy criterion in frequency-wavenumber inversion: Geophysics, **56**, 2027-2035.

Billette, F., and G. Lambare, 1998, Velocity macro-model estimation from seismic reflection data by stereotomography: Geophysical Journal International, **135**, 671-690.

Brenders, A., and R. Pratt, 2007, Full waveform tomography for lithospheric imaging: Results from a blind test in a realistic crustal model: Geophysical Journal International, **168**, 133-151.

Bunks, C., F. Saleck, S. Zaleski, and G. Chavent, 1995, Multiscale seismic waveform inversion: Geophysics, **60**, 1457-1473.

Crase, E., A. Pica, M. Noble, J. McDonald, and A. Tarantola, 1990, Robust elastic nonlinear waveform inversion: Application to real data: Geophysics, **55**, 527-538.

Guitton, A., and W. Symes, 2003, Robust inversion of seismic data using the Huber norm: Geophysics, **68**, 1310-1319.

Ha, W., W. Chung, E. Park, and C. Shin, 2012a, 2-D acoustic Laplace-domain waveform inversion of marine field data: Geophysical Journal International, **190**, 421-428.

Ha, W., W. Chung, and C. Shin, 2012b, Pseudo-Hessian matrix for the logarithmic objective function in full waveform inversion: Journal of Seismic Exploration, **21**, 201-214.

Ha, W., and C. Shin, 2013, Why do Laplace-domain waveform inversions yield long-wavelength results?: Geophysics, **78**, R167-R173.

Ha, W., J. Yoo, and C. Shin, 2011, Efficient velocity estimation in the Laplace domain using gain control: SEG Expanded Abstracts, 4071-4076.

Koo, N.-H., C. Shin, D.-J. Min, K.-P. Park, and H.-Y. Lee, 2011, Source estimation and direct wave reconstruction in Laplace-domain waveform inversion for deep-sea seismic data: Geophysical Journal International, **187**, 861-870.

Kreyszig, E., 2011, Advanced Engineering Mathematics (10 ed.): John Wiley & Sons, Inc.

Kwak, S., H. Jun, W. Ha, and C. Shin, 2013, Temporal windowing and inverse transform of the wavefield in the Laplace-Fourier domain: Geophysics, **78**, R207-R222.

Officer, C. B., 1958, Introduction to the theory of sound transmission : McGraw-Hill.

Operto, S., J. Virieux, J. X. Dessa, and G. Pascal, 2006, Crustal seismic imaging from multifold ocean bottom seismometer data by frequency domain full waveform tomography: Application to the eastern Nankai trough: Journal of Geophysical Research, **111**, B09306.

Park, E., W. Ha, W. Chung, C. Shin, and D.-J. Min, 2013, 2D Laplace-Domain Waveform Inversion of Field Data Using a Power Objective Function: Pure and Applied Geophysics,

Pratt, R., C. Shin, and G. Hicks, 1998, Gauss-Newton and full Newton methods in frequency-space seismic waveform inversion: Geophysical Journal International, **133**, 341-362.

Pratt, R. G., Z. M. Song, P. Williamson, and M. Warner, 1996, Two-dimensional velocity models from wide-angle seismic data by wavefield inversion: Geophysical Journal,

**124**, 323-340.

Press, W., S. Teukolsky, W. Vetterling, and B. Flannery, 1992, Numerical Recipes in Fortran 77 (2nd ed.): Cambridge University Press.

Ravaut, C., S. Operto, L. Improta, J. Virieux, A. Herrero, and P. Dell'Aversana, 2004, Multiscale imaging of complex structures from multifold wide-aperture seismic data by frequency-domain full-waveform tomography: application to a thrust belt: Geophysical Journal International, **159**, 1032-1056.

Shin, C., and Y. Cha, 2008, Waveform inversion in the Laplace domain: Geophysical Journal International, **173**, 922-931.

Symes, W., 2008, Migration velocity analysis and waveform inversion: Geophysical Prospecting, **56**, 765-790.

Tarantola, A., 1984, Inversion of seismic-reflection data in the acoustic approximation: Geophysics, **49**, 1259-1266.

Virieux, J., and S. Operto, 2009, An overview of full-waveform inversion in exploration geophysics: Geophysics, **74**, WCC1-WCC26.

Zelt, C., G. Pratt, A. Brenders, S. Hanson-Hedgecock, and J. A. Hole, 2005, Advancements in long-offset seismic imaging: A blind test of traveltime and waveform tomography: AGU Spring meeting, S52A-S504.

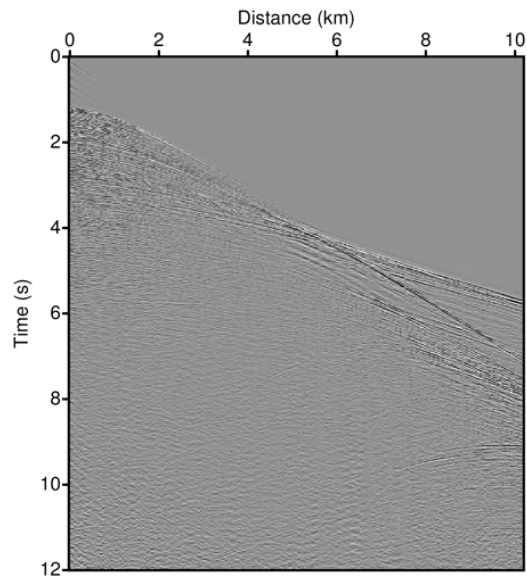 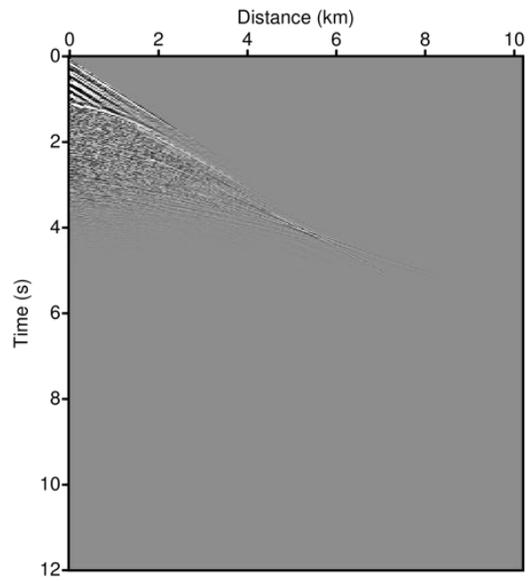

(a) (b)

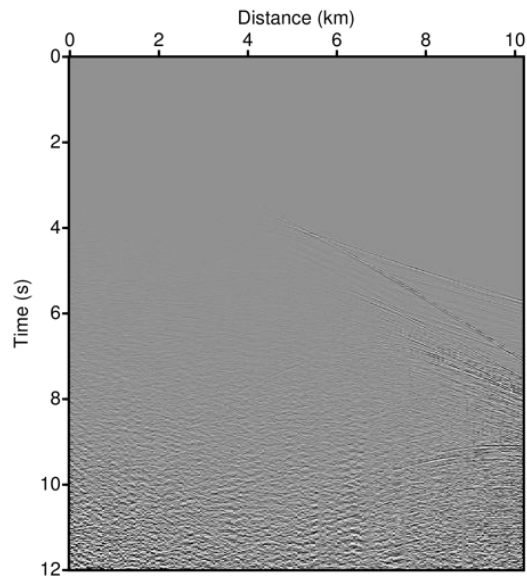 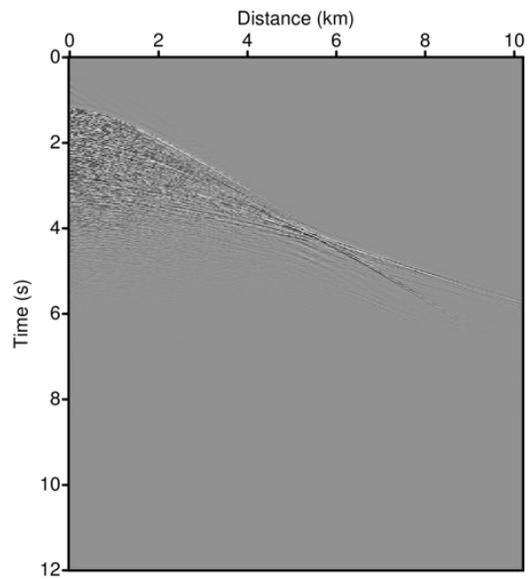

(c) (d)

Figure 1. (a) A shot gather from the Gulf of Mexico data and the shot gathers after applying a function of (b) $e^{-2t}$, (c) $t^4$, and (d) $t^4 e^{-2t}$.

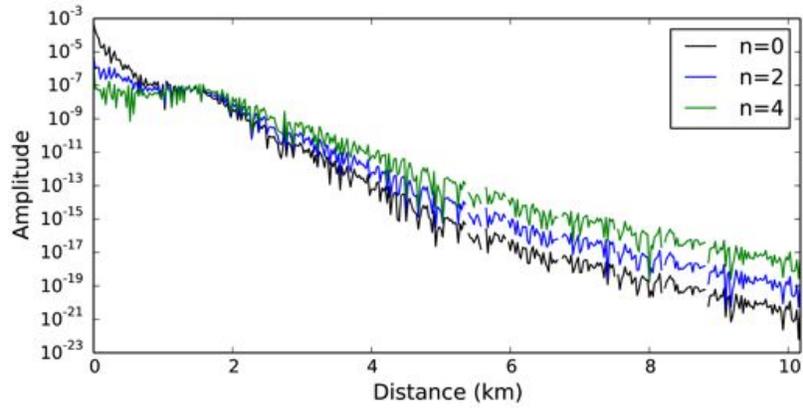

Figure 2. Laplace-transformed shot gathers with the power of 0, 2, and 4 when the damping constant is 7 $s^{-1}$.

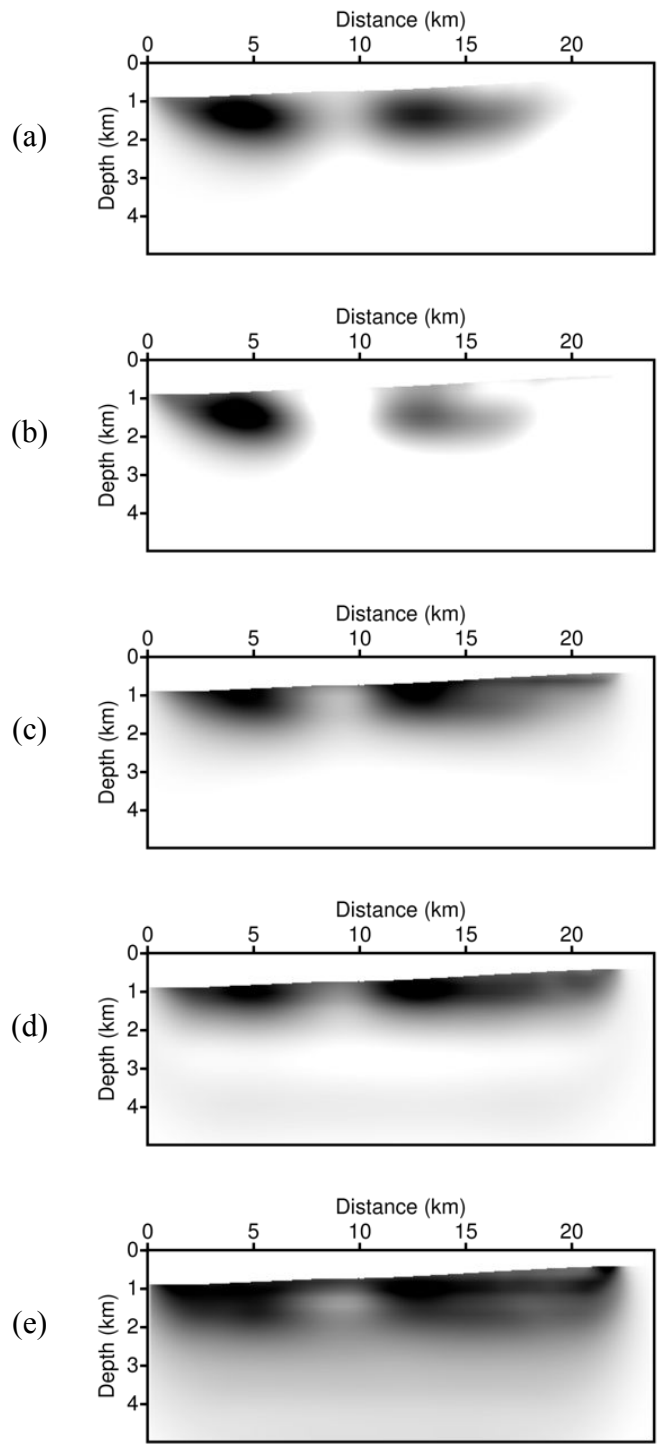

Figure 3. The gradient directions obtained using equation 6 with (a) n=0, (b) n=1, (c) n=2, (d) n=3, and (e) n=4.

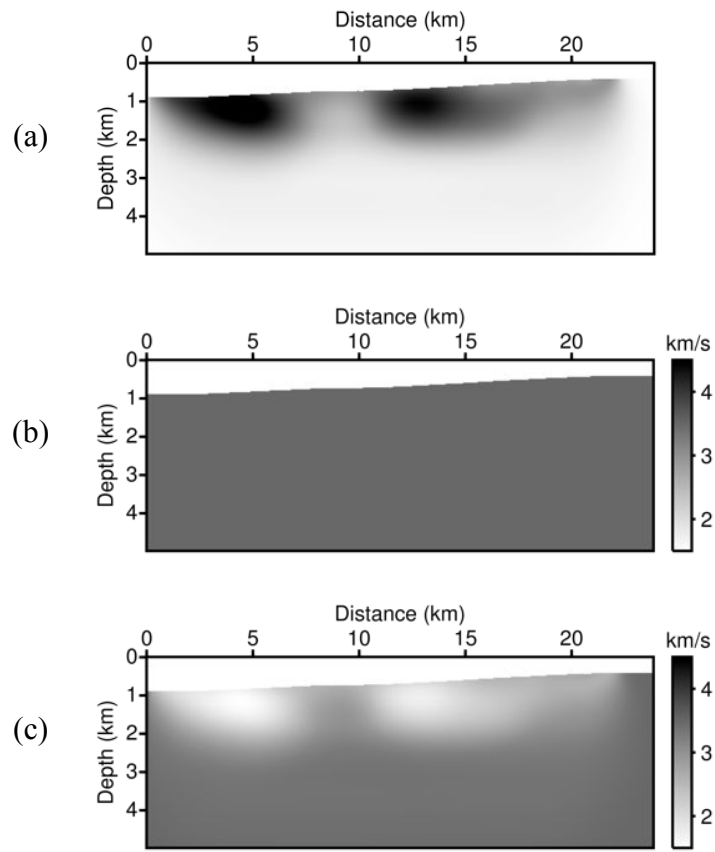

Figure 4. (a) The final gradient direction, (b) the homogeneous initial velocity model with the velocity of 3.5 km/s at the water bottom, and (c) the velocity model updated once from the homogeneous model using the final gradient.

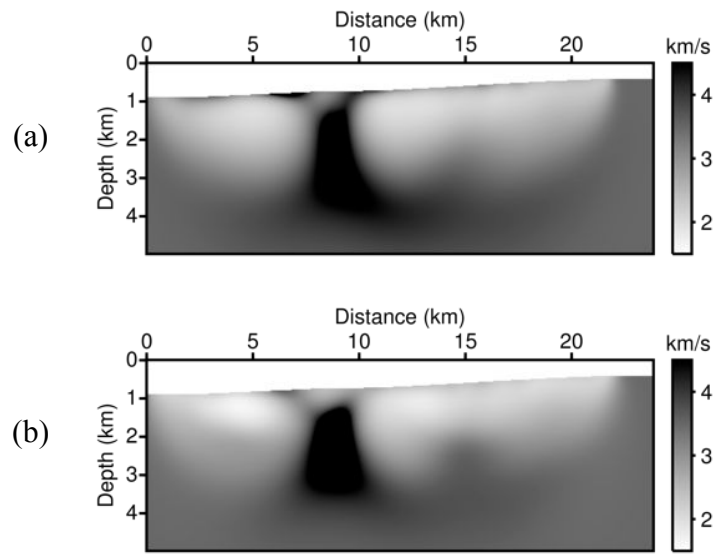

Figure 5. Inversion results obtained after 80 iterations started from (a) the homogeneous initial model (Figure 4b), and (b) the updated model (Figure 4c).

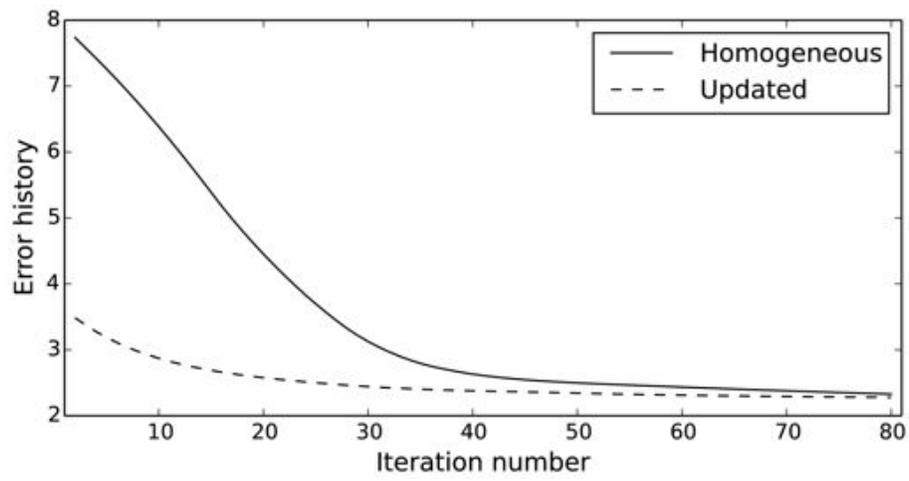

Figure 6. The error histories of the two inversion examples.

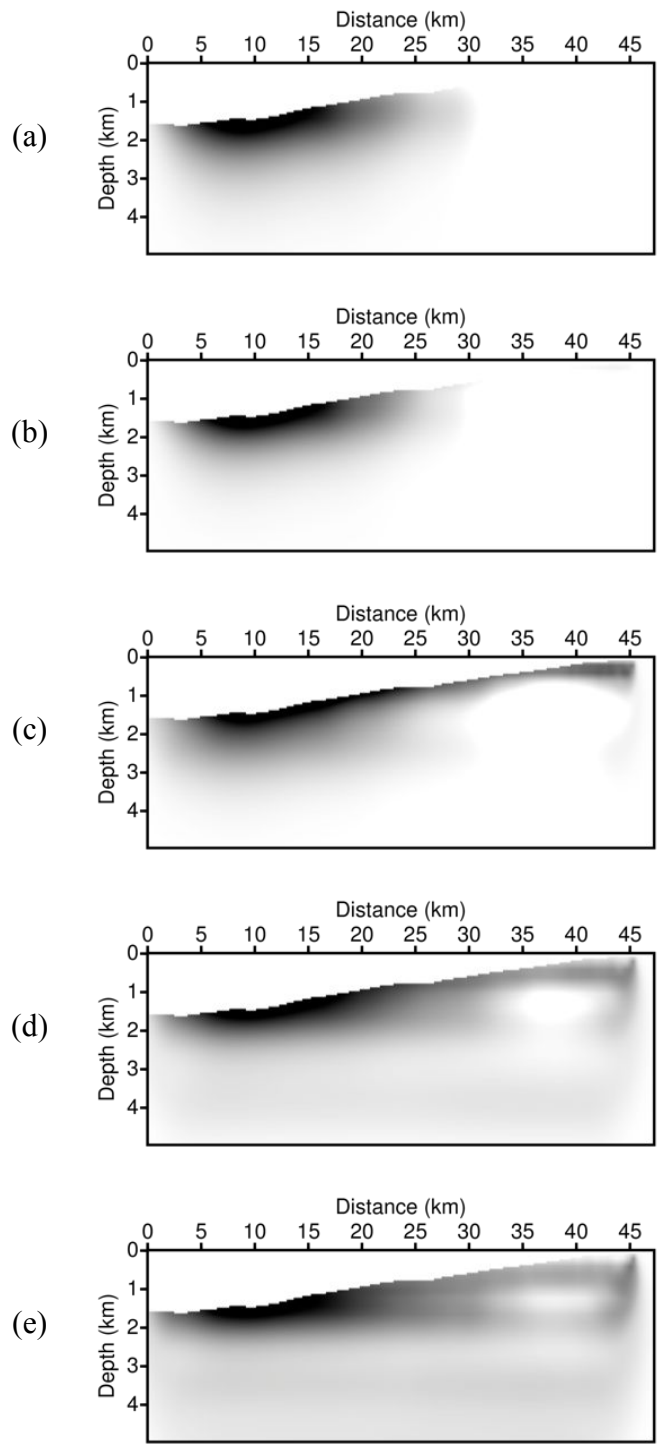

Figure 7. The gradient directions obtained using equation 6 with (a) n=0, (b) n=1, (c) n=2, (d) n=3, and (e) n=4.

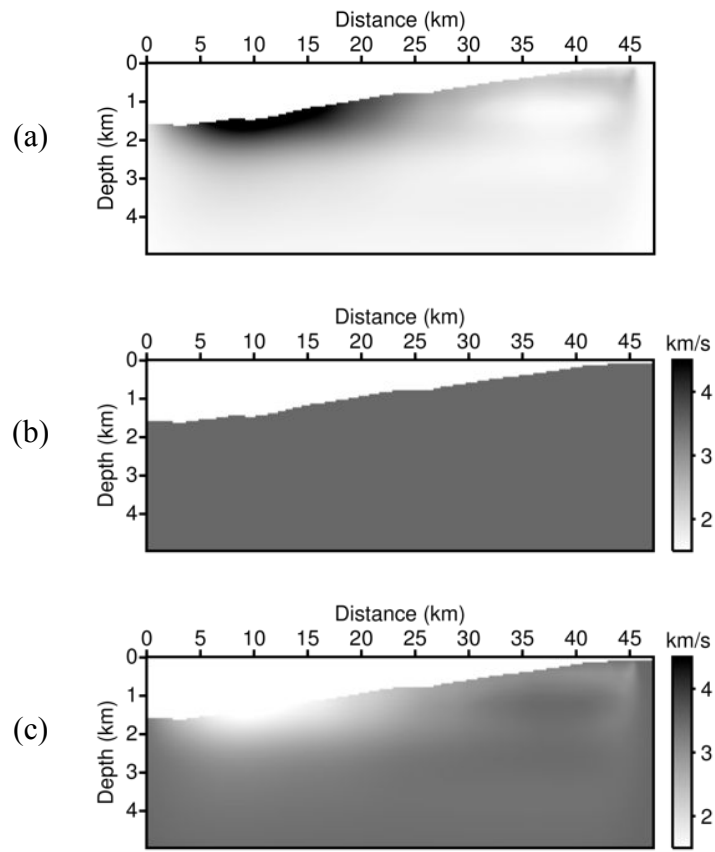

Figure 8. (a) The final gradient direction, (b) the homogeneous initial velocity model with the velocity of 3.5 km/s at the water bottom, and (c) the velocity model updated once from the homogeneous model using the final gradient.

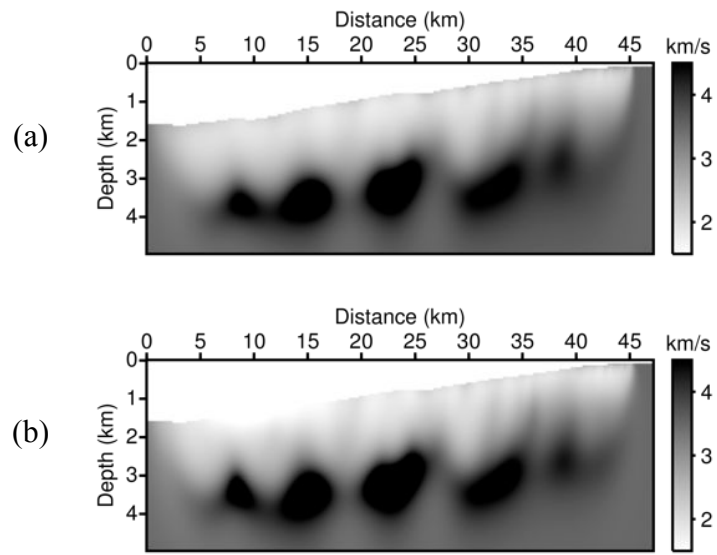

Figure 9. Inversion results obtained after 80 iterations started from (a) the homogeneous initial model (Figure 8b), and (b) the updated model (Figure 8c).

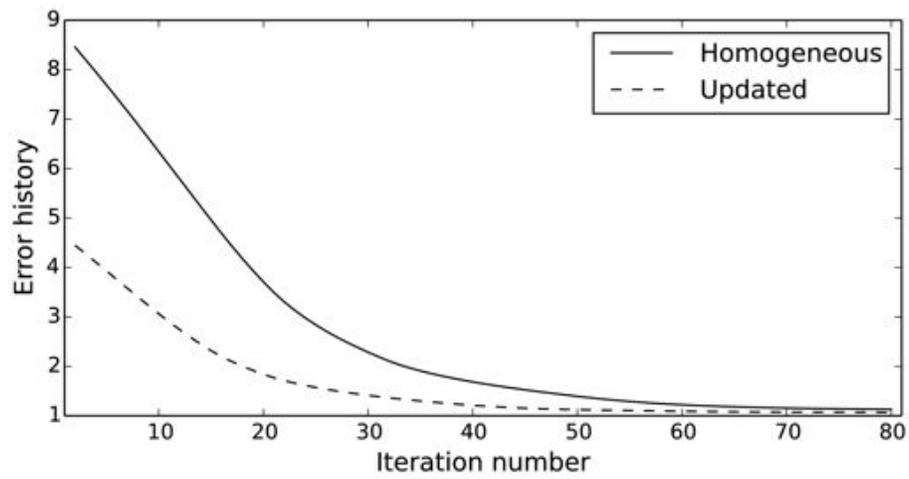

Figure 10. The error histories of the two inversion examples.

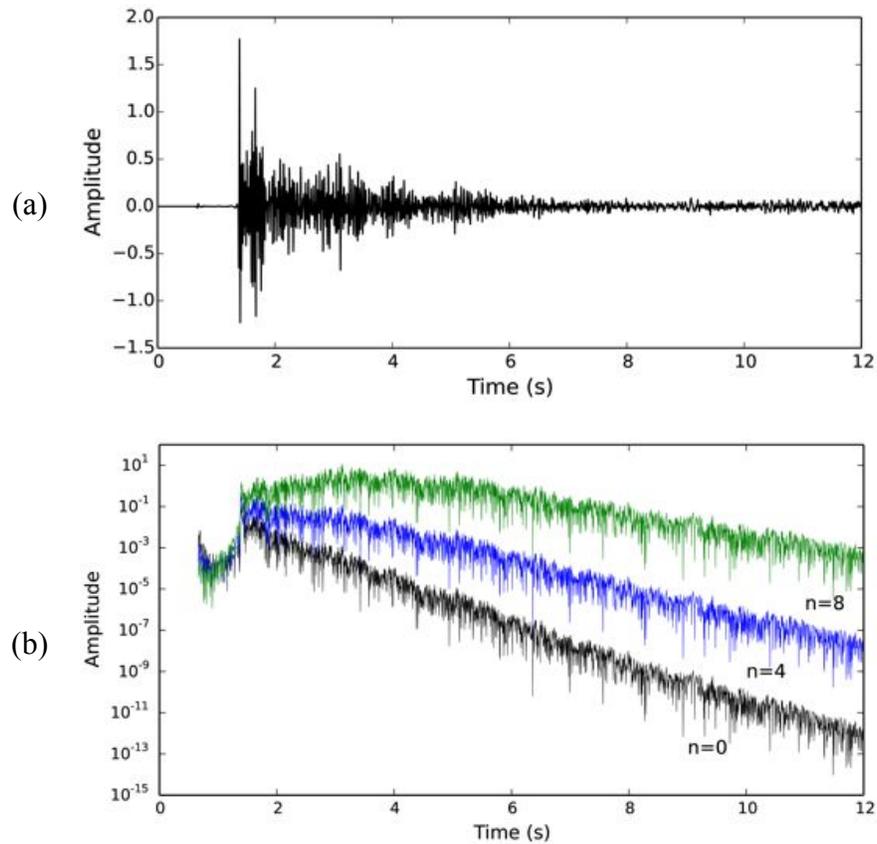

Figure 11. (a) A trace from the shot gather shown in figure 1 with the offset of 1,021 m and (b) the trace damped with the damping constant of 2 $s^{-1}$ for different powers of the gain function. The absolute values of the traces are shown in the logarithmic scale.